\documentclass[11pt]{amsart}

\usepackage {amsfonts,amsthm}
\usepackage[english]{babel}

\theoremstyle{plain}
\newtheorem{Th}{Theorem}
\newtheorem{Le}{Lemma}

\def\R{{\mathbb{R}}}
\def\Z{{\mathbb{Z}}}
\def\N{{\mathbb{N}}}
\def\s{\sigma}
\def\d{\delta}
\def\f{\varphi}

\def\o{\omega}
\def\e{{\varepsilon}}
\def\Lin{{\rm Lin}}

\def\supp{{\rm supp}}
\def\Fr{{\it Fr}}

\begin{document}
\title
{Bohr and Besicovitch almost periodic discrete sets\\ and quasicrystals}

\author{Favorov S.}

\address{Mathematical School, Kharkov National University, Swobody sq.4,
Kharkov, 61077 Ukraine}

 \email{sfavorov@gmail.com}

\date{}

\begin{abstract}
 A discrete set $A$ in the Euclidian space is  almost periodic if  the
 measure with the unite masses at points of the set is almost periodic in the weak sense.
We investigate properties of such sets in the case when $A-A$ is discrete. In particular,
if $A$ is a Bohr almost periodic set, we prove that $A$ is a union of a finite
number of translates of a certain full--rank lattice. If $A$ is a Besicovitch almost
periodic set, then there exists a full-rank lattice such that in most cases a nonempty
intersection of its translate with $A$ is large.
\end{abstract}



\maketitle

The notion of an almost periodic discrete set in the complex plane is well known in
the theory of almost periodic holomorphic and meromorphic functions
(cf.~\cite{L},\cite{T}, \cite{FRR}, \cite{FP}, \cite{F}). Almost periodic discrete
sets in the $p$-dimensional Euclidian space appear later in the mathematical theory of
quasicrystals (cf.~\cite{La},\cite{M}), where various notions of almost periodicity
were used.

In this connection, the following question was raised in \cite{La} (Problem 4.4): whether any
Bohr almost periodic discrete set is a finite union of translations of a full--rank
lattice in $\R^p$. In \cite{FK} and \cite{FK1} we showed that every almost periodic
perturbation of a full--rank lattice in $\R^p$ is a Bohr almost periodic set. Hence,
there exists a wide class of such sets.

Next, most of the models of quasicrystals are sets of a finite type. In the
present article we prove that every Bohr almost periodic discrete set of a finite type is
just a finite union of translations of a full--rank lattice. We also prove that every
Besicovitch almost periodic discrete set of a finite type has, in a sense,
asymptotically the same form.
\medskip

Let us recall some known definitions (see, for example, \cite{C}, \cite{R}).

\medskip
A continuous function $f(x)$ in $\R^p$ is {\it almost periodic in the sense of Bohr}
if for any $\e>0$ the set of {\it $\e$-almost periods of $f$}
 $$
  E_\e=\{\tau\in\R^p:\,\sup_{x\in\R^p}|f(x+\tau)-f(x)|<\e\}
  $$
  is a relatively dense set in $\R^p$. The latter means that there is
  $R=R(\e)<\infty$ such that every ball of radius $R$ contains at least one $\e$-almost period of
  $f$.

The following definition is an evident generalization of the one given by Besicovitch
 in \cite{B} for the
case when $p=1$:
an integrable function $F(x)$ in $\R^p$ is {\it almost periodic in the sense of
Besicovitch, or $B^1$-almost periodic}, if for any $\e>0$ there is a Bohr almost
periodic function $f_\e(x)$ in $\R^p$ such that
 $$
\limsup_{R\to\infty}\frac{1}{\o_pR^p}\int_{|t|<R}|F(x)-f_\e(x)|~dx<\e.
 $$
Here $\o_p$ is the volume of the unit ball in $\R^p$. Next, $T\in\R^p$ is an
$\e$-almost period of  $F$, if
 $$
\limsup_{R\to\infty}\frac{1}{\o_pR^p}\int_{|t|<R}|F(x+T)-F(x)|~dx<\e.
 $$
It follows from the definition of $B^1$-almost periodicity that the set $E_\e$ of
$\e$-almost periods of $F$ is a relatively dense set in $\R^p$ for each $\e>0$.

\medskip
We will say that a set $A\subset\R^p$ is {\it discrete}, if it has no finite limit
points, and {\it uniformly discrete}, if there is $r>0$ such that any ball of radius $r$
contains at most one point of $A$.  A discrete set $A\subset\R^p$ is {\it Bohr
(Besicovitch) almost periodic}, if for every continuous function $\f$ in $\R^p$ with a
compact support the sum $\sum_{a\in A}\f(x-a)$ is a Bohr (Besicovitch) almost
periodic function, i.e., the union of the unit masses at the points of $A$ is a Bohr
(Besicovitch) almost periodic measure in the weak sense.

  \smallskip
 There is a geometric criterion for a discrete set to be
Bohr almost periodic.
\begin{Th}[\cite{FK}, Theorem 11]\label{FK}
 A discrete set $A=(a_n)_{n\in\N}\subset\R^p$ is Bohr almost  periodic if and only if for each
 $\e>0$ the set of $\e$-almost  periods of $A$
  $$
\{\tau\in\R^p:\,\exists \hbox{\ a bijection}\quad \s:\,\N\to\N\quad\hbox{ such that
}\quad |a_n+\tau-a_{\sigma(n)}|<\e\quad \forall\,n\in\N\}
  $$
 is relatively dense in $\R^p$.
\end{Th}

It follows easily from this criterion  that the number of elements of $A$ in every ball
of radius $1$ is uniformly bounded.

 Note that for a uniformly discrete set $A$ and
sufficiently small $\e>0$, a vector $\tau$ is an $\e$-almost  period of $A$ if and only
if for any $a\in A$ there exists $a'\in A$ such that $|a+\tau-a'|<\e$.

\medskip
A set $A$ is {\it a Delone set,} if it is uniformly discrete and relatively dense in
$\R^p$. For sets $A,\, C\subset\R^p$ put $A+C=\{a+c:\,a\in A,\,c\in C\}$, $A-C$ is defined in a similar way. Also, put $B(x,R)=\{y\in\R^p:|y-x|<R\}$. Following
\cite{La}, we will say that a discrete set $A\subset\R^p$ is  {\it a set of finite
type}, if the set $A-A$ is discrete. A set $A$ is {\it a Meyer set}, if the set $A-A$
is a Delone set. Clearly, any Mayer set is a Delone set. A set $A\subset\R^p$ is an
{\it ideal crystal}, if $A$ consists of a finite number of translates of a full--rank
lattice $L$, that is, $A=L+F$, where $F$ is a finite set and $L$ is an additive
discrete subgroup of $\R^p$ such that $\Lin_\R L=\R^p$.

We shall prove the following theorem.
\begin{Th}
  If a set $A$ of a finite type is Bohr almost periodic, then it is an ideal crystal.
\end{Th}

{\bf Proof}. Every ball of radius $R=R(1)$ contains at least one $1$-period of $A$,
hence each ball of radius $R+1$ contains at least one point $a\in A$. Since the set
$A-A$ is discrete, we see that there is $\e>0$ such that $\e<\min\{1;|(a-b)-(c-d)|\}$
whenever $a,\,b,\,c,\,d\in A$ and $|a-b|<2R+4,\,|c-d|<2R+4,\,a-b\neq c-d$.  In
particular, $\e<|a-b|$ whenever $a,\,b\in A$ and $a\neq b$.

Fix $a\in A$. Let $\tau\in\R^p$ be an arbitrary $(\e/2)$-almost period of $A$. Taking
into account our choice of $\e$, we see that there is a unique $c\in A$ that satisfies
the inequality $|a+\tau-c|<\e/2$. Clearly, $T=c-a$ is an $\e$-almost period of $A$.
Let us show that $T$ is actually a period of $A$.

Suppose that $b\in A$ such that $b\neq a$ and $|a-b|<2R+3$. Since $T$ is an
$\e$-almost period of $A$, there exists a point $d\in A$ such that
$|b+T-d|=|(a-b)-(c-d)|<\e$. Since $|c-d|\le|a-b|+|b+T-d|<2R+4$, we obtain $a-b=c-d$
and $d=b+T$. We repeat these arguments for all $b\in A$ such that $|b-a|<2R+3$ and, after that, for all $b'\in A$ such that $|b'-b|<2R+3$. After countable number of steps we obtain
that $a+T\in A$ for all $a\in A_1\subset A$. If $A\setminus A_1\neq\emptyset$, then
$R_1=\inf\{|a-b|:\,a\in A_1,\,b\in A\setminus A_1\}\ge 2R+3$. Take $a'\in A_1$ and
$b'\in A\setminus A_1$ such that $|a'-b'|<R_1+1$. This implies that
$B(\frac{a'+b'}{2},R+1)\cap A_1=\emptyset$ and $B(\frac{a'+b'}{2},R+1)\cap(A\setminus
A_1)=\emptyset$, which is impossible. Hence, $A=A_1$ and $T$ is a period of $A$.

 Next, take $(\e/2)$-almost periods $\tau_j$ from the set
 $$
 \left\{x\in\R^p:\,|x|>4p,\,|x-\langle x,e_j\rangle|<\frac{|\langle x,e_j\rangle|}{4p}\right\},
 \quad j=1,\dots,p,
 $$
where $e_j,\,j=1,\dots,p,$ is a basis in $\R^p$. There are periods
$T_j$ such that $|T_j-\tau_j|<1/2$. Since $|\langle\tau_j,e_j\rangle|>2p$, we get
 \begin{equation}
 \frac{\max_{k\neq j}|\langle T_j,e_k\rangle|}{|\langle T_j,e_j\rangle|}\le
 \frac{|T_j-\langle T_j,e_j\rangle|}{|\langle T_j,e_j\rangle|}<
 \frac{|\tau_j-\langle\tau_j,e_j\rangle|+1}{|\langle\tau_j,e_j\rangle|-1/2}<\frac{1}{p},
 \quad j=1,\dots,p.
 \end{equation}
 Hence, the determinant of the matrix
 $(\langle T_j,e_k\rangle)_{j,k=1}^p$ does not vanish, and the vectors
$T_1,\dots,T_p$ are linearly independent. Consequently, the set $L=\{n_1T_1+\dots
+n_pT_p:\,n_1,\dots,n_p\in\Z\}$ is a full--rank lattice. Next, the set $F=\{a\in A:
|a|<|T_1|+\dots+|T_p|\}$ is finite. All vectors $t\in L$ are periods of $A$, hence,
$L+F\subset A$. On the other hand, for each $a\in A$ there is $t\in L$ such that
$|a-t|<|T_1|+\dots+|T_p|$, hence, $a-t\in F$. The theorem is proved.

\medskip
 Let $A$ be a discrete set and let $L$ be a full--rank lattice. For $a, b\in A$ put $a\sim b$ if
 $a-b\in L$. Then there is a unique at most countable decomposition $A=\cup_j A_j^L$  into mutually
 disjoint equivalence classes $A_j^L$. Next, by  $\# E$ denote the number of elements in the finite set $E$.

\medskip
\begin{Th}
  Let a Besicovitch almost periodic discrete set $A\subset\R^p,\,p>1,$ be a Meyer set. Then for any
  $N<\infty$ and $\eta>0$ there is a full--rank lattice $L$ such that
 $$
\sum_{j: \#(A_j^L\cap B(0,R))>N}\#(A_j^L\cap B(0,R))\ge(1-\eta)\#(A\cap B(0,R)),\qquad
R>R(N,\eta).
 $$
  \end{Th}
Our proof is based on the following lemmas.

\begin{Le}
 Let a Delone set $A\subset\R^p$ be Besicovitch almost periodic. Then for any $\e>0$ and
 $\d>0$
 there is a relatively dense set $E=E(\e,\d)\subset\R^p$ such that for any $\tau\in E$
 \begin{equation}\label{1}
\forall a\in A\setminus\tilde A\quad\exists a'\in A:\quad|a+\tau-a'|<\e,
 \end{equation}
where $\tilde A=\tilde A(\tau)\subset A$ such that for sufficiently large $R$
 \begin{equation}\label{2}
\#(\tilde A\cap B(0,R))<\d\#(A\cap B(0,R)).
 \end{equation}
 \end{Le}

{\bf Proof}. The set $A$ is relatively dense, therefore $A\cap B(x,R_0)\neq\emptyset$
for some $R_0<\infty$ and all $x\in\R^p$. Hence there exists $\kappa>0$ such that for
all sufficiently large $R$
 \begin{equation}\label{0}
\# A\cap B(0,R)>\kappa R^p.
 \end{equation}
We may suppose that $\e<\frac{1}{2}\inf\{|a-b|:\,a,b\in A,\,a\neq b\}$. Let $\f$ be a
$C^\infty$ function in $\R^p$ such that $0\le\f(x)\le1$, $\f(x)=1$ for $|x|<1/4$ and
$\f(x)=0$ for $|x|>1/2$. Clearly, $\int\f(x)~dx>\o_p4^{-p}$. The function
$\psi(x)=\sum_{a\in A}\f(\frac{x-a}{\e})$ is  a Besicovitch almost periodic function.
Clearly, if $\psi(x)>0$ for $x\in\R^p$, then $x\in B(a,\e/2)$ for some $a\in A$. Let
$E$ be the set of $\d\kappa(\e/4)^p$-almost periods of $\psi$ and $-\tau\in E$. Put
$\tilde A=\{a\in A:\,B(a+\tau,\e)\cap A=\emptyset\}$. Since $\tilde A\subset
A\setminus\supp~\psi(x-\tau)$, we get that for $R>R(\tau)$
 $$
\o_p(\e/4)^p\#(\tilde A\cap B(0,R))<\sum_{a\in\tilde A\cap
B(0,R)}\int_{B(a,\e/2)}\f\left(\frac{x-a}{\e}\right)~dx
 $$
 $$
\le \int_{B(0,R+\e)\setminus\supp~\psi(x-\tau)}\psi(x)~dx\le
\int_{B(0,R+\e)}|\psi(x-\tau)-\psi(x)|~dx<(\e/4)^p\d\kappa\o_pR^p.
 $$

Hence,
 $$
\#(\tilde A\cap B(0,R))<\d\kappa R^p.
 $$
The assertion of the lemma follows from (\ref{0}).

\begin{Le}
Suppose that $E$ is a finite subset of a full--rank lattice
$L=\Lin_\Z\{T_1,\dots,T_p\}\subset\R^p,\,p>1$. Put
 $$
 \Fr(E)=\{b\in E:\,\hbox{at least one of the points $b\pm T_1,\dots,b\pm T_p$ does not belong to
  $E$}\}.
 $$
Then $(\# E)^{\frac{p-1}{p}}\le (p-1)\# \Fr(E)$.
\end{Le}
{\bf Proof}. Without loss of generality we may suppose $L=\Z^p$. For any
$F\subset\Z^p$ denote by $P_j(F)$ the projection of $F$ on the hyperplane $x_j=0$.
Note that $\# \Fr(E)\ge\# P_j(E)$ for all $j\in\overline{1,p}$.

Put
 $$
 E_1=\{x\in E:\,\# E\cap P_1^{-I}P_1(\{x\})\le(\# E)^{1/p}\},
  $$
and, for each  $j\in\overline{2,p-1}$ recurrently
 $$
 E_j=\{x\in(E\setminus\cup_{i=1}^{j-1}E_i):\,\#(P_j^{-I}P_j(\{x\})\cap(E\setminus\cup_{i=1}^{j-1}E_i))
 \le(\# E)^{1/p}\}.
  $$
Clearly, $\# E_j\le\#P_j(E_j)(\# E)^{1/p}$. If $\# E_j\ge(\# E)/(p-1)$ for some
$j\in\overline{1,p-1}$, we get
 $$
\#\Fr(E)\ge\# P_j(E)\ge\# P_j(E_j)\ge (p-1)^{-1}(\# E)^{1-1/p}.
 $$
If $\# E_j<(p-1)^{-1}(\# E)$ for all $j\in\overline{1,p-1}$, then there exists $x'\in
E\setminus\cup_{j=1}^{p-1}E_j$. Let $x'_p$ be the last coordinate of $x'$. Since
$x'\not\in E_{p-1}$, we see that the set $P_{p-1}^{-I}P_{p-1}(\{x'\})\cap(E\setminus
\cup_{i=1}^{p-2}E_i)$ contains at least $(\# E)^{1/p}$ points. Since every point $x''$
from the latter set  does not belong to $E_{p-2}$, we get that all the sets
$P_{p-2}^{-I}P_{p-2}(\{x''\}) \cap(E\setminus\cup_{i=1}^{p-3}E_i)$ contain at least
$(\# E)^{1/p}$ points. Continuing this line of reasoning, we see that there exist at
least $(\# E)^{\frac{p-1}{p}}$ different points in $E$ with the last coordinate
$x'_p$. Hence, these points have  distinct projections on the hyperplane $x_p=0$.
Consequently, in this case we get $\#\Fr(E)\ge\# P_p(E)\ge(\# E)^{\frac{p-1}{p}}$.
Lemma follows.
\medskip

{\bf Proof of the theorem}. Note that $A$ is a Delone set. Applying Lemma 1 with
$\e<\frac{1}{2}\inf\{|x-y|:\,x,y\in A-A,\,x\neq y\}$ and
$\d=2^{-1}p^{-2}N^{-1/p}\eta$, we find $\tau\in\R^p$ and $A(\tau)\subset A$ such that
(\ref{1}) and (\ref{2}) hold. Replace $\tau$ by $T=a'-a$. Taking into account the
inequality $|a+\tau-a'|<\e$ and the bound on $\e$, we get
 $$
 \forall a\in A\setminus\tilde A(\tau)\quad\exists a'\in A:\quad a'=a+T,
 $$
 Moreover,  arguing as at the end of the proof of
Theorem 2, we can take $\e/2$-almost periods $\tau_1,\dots,\tau_p$ that are linearly
independent over $\R$, vectors $T_1,\dots,T_p$, and $A^*=\cup_j\tilde
A(\tau_j)\cup\cup_j\tilde A(-\tau_j)$ such that
 $$
 \forall a\in A\setminus A^*\quad\exists a'_j, a''_j\in A:\quad
 a'_j=a+T_j,\quad a''_j=a-T_j,\quad j\in\overline{1,p},
 $$
 where the set $A^*$  satisfies the following bound for all sufficiently large $R$
 \begin{equation}\label{a}
\# A^*\cap B(0,R)<2p\d(\# A\cap B(0,R))=p^{-1}N^{-1/p}\eta(\# A\cap B(0,R)).
 \end{equation}
 Put $L=\{n_1T_1+\dots +n_pT_p:\,n_1,\dots,n_p\in\Z\}$ and $J=\{j:\,\#(A_j^L\cap
B(0,R))\le N\}$. Taking into acount Lemma 2 and the definition of $J$, we get
 \begin{equation}\label{b}
 \sum_{j\in J}\#(A_j^L\cap B(0,R))\le(p-1)\sum_{j\in J}\#\Fr(A_j^L\cap B(0,R))N^{1/p}.
 \end{equation}
 Every  point $a\in\Fr(A_j^L\cap B(0,R))$ such that $|a|+\max_j|T_j|<R$ belongs to  $A^*$, therefore,
 \begin{equation}\label{c}
\sum_{j\in J}\#\Fr(A_j^L\cap B(0,R))\le\#[A^*\cap B(0,R)]+\# A\cap[B(0,R)\setminus
B(0,R-\max_j|T_j|)].
 \end{equation}
 Since $A$ is a uniformly discrete and relatively dense set, we get
 \begin{equation}\label{d}
 \# A\cap[B(0,R)\setminus B(0,R-\max_j|T_j|)]=o(\# A\cap B(0,R))\quad\hbox{as}\quad
 R\to\infty.
 \end{equation}
 It follows from (\ref{a}---\ref{d})
 that
 $$
 \sum_{j\in J}\#(A_j^L\cap B(0,R))\le\eta\#(A\cap B(0,R)),\quad R>R(N,\eta).
 $$
The last inequality yields the assertion of the theorem.

\medskip
We see that different notions of almost periodicity lead to different results. In
particular, it looks like the almost periodicity in the sense of Bohr for discrete sets is too restrictive and because of that it probably will have limited applications
in mathematical theory of quasicrystals.

\end{document}